\documentclass[12pt,leqno,draft]{article}

\newtheorem{theorem}{Theorem}
\newtheorem{lemma}[theorem]{Lemma}
\newtheorem{proposition}[theorem]{Proposition}
\newtheorem{definition}[theorem]{Definition}
\newtheorem{corollary}[theorem]{Corollary}

\newcommand{\begintheorem}{\addtocounter{equation}{1}\begin{theorem}}
\newcommand{\beginlemma}{\addtocounter{equation}{1}\begin{lemma}}
\newcommand{\beginproposition}{\addtocounter{equation}{1}\begin{proposition}}
\newcommand{\begindefinition}{\addtocounter{equation}{1}\begin{definition}}
\newcommand{\begincorollary}{\addtocounter{equation}{1}\begin{corollary}}

\begin{document}

\title{Some topics in complex and harmonic analysis, 2}

\author{Stephen William Semmes	\\
	Rice University		\\
	Houston, Texas}

\date{}

\maketitle

	Let $h(x)$ be a complex-valued bounded continuous function on
${\bf R}^n$.  For each positive real number $t$, define $A_t(x)$ to be
$\exp ( - 2 \pi t \, |x|)$, $x \in {\bf R}^n$.  The corresponding Abel
mean of $h$ is defined to be
\begin{equation}
	\int_{{\bf R}^n} h(x) \, A_t(x) \, dx.
\end{equation}
If the limit of this as $t \to 0$ exists, then we say that $h(x)$
is Abel summable with the limit as the sum.  If $h(x)$ is integrable
on ${\bf R}^n$, then $h$ is Abel summable with the usual integral of
$h$ being the Abel sum.

	One can show that the function $A_t(x)$ can be expressed as an
average of Gaussians.  This is a quite useful fact, and basically it
is a way of saying that Gaussians are more concentrated.  For instance,
one can use this to show that if $h$ is Gauss summable, then $h$ is
Abel summable, and with the same sum.  One can also use this to compute
the Fourier transform of $A_t$, which is the Poisson kernel
\begin{equation}
	P_t(y) = c_n \frac{t}{(|y|^2 + t^2)^{(n+1)/2}},
\end{equation}
where $c_n = \Gamma((n+1)/2) \, \pi^{-(n+1)/2}$.

	Now suppose that $f(x)$ is a continuous integrable function
on ${\bf R}^n$, so that the Fourier transform of $f$ is defined by
\begin{equation}
	\widehat{f}(\xi) 
	   = \int_{{\bf R}^n} f(x) \, \exp (-2 \pi i \, \xi \cdot x) \, dx,
\end{equation}
with $\xi \cdot x = \sum_{j=1}^n \xi_j \, x_j$.  The Fourier transform
of $f$ is a bounded continuous function, and even uniformly
continuous.  It may not be integrable, and we can always consider the
Abel sums for the inverse Fourier transform,
\begin{equation}
\label{abel sum for inverse fourier transform of f}
  	\int_{{\bf R}^n} \widehat{f}(\xi) \, 
		A_t(\xi) \, \exp (2 \pi i \, x \cdot \xi) \, d\xi.
\end{equation}
One can show that this is equal to $(P_t * f)(x)$, where the
convolution $(P_t *f)(x)$ is defined as usual as the integral of
$P_t(y)$ times $f(x - y)$ with respect to $y$.  Also, $(P_t * f)(x)
\to f(x)$ as $t \to 0$, uniformly on compact subsets of ${\bf R}^n$.

	Now suppose that $n = 1$, and that $f(x)$ is a continuous
integrable function on the real line such that the Fourier transform
of $f$ is supported in $[0, \infty)$.  If $z = x + i \, y$ is a
complex number, $x, y \in {\bf R}$, which lies in the upper
half-plane, so that $y > 0$, then consider the integral
\begin{equation}
	\int_0^\infty \widehat{f}(\xi) \, \exp (2 \pi i \, z \, \xi) \, d\xi.
\end{equation}
This is in fact the same as the Abel sum (\ref{abel sum for inverse
fourier transform of f}) for the inverse Fourier transform of $f$ with
$y = t$.

	Writing the integral in this manner suggests strongly that
this defines a holomorphic function of $z$ in the upper half-plane,
which can be verified through standard arguments.

	Let us denote this function $f(z)$, for $z$ in the upper
half-plane.  Of course we also have $f(z)$ defined when $y = 0$, which
is to say that $z = x \in {\bf R}$.  In other words, $f(z)$ is defined
for all complex numbers $z = x + i \, y$, $x, y \in {\bf R}$, such
that $y \ge 0$.  Using the continuity of $f$ on ${\bf R}$ and the fact
that $f(z) = (P_y * f)(x)$ when $y > 0$ one can show that $f(z)$ is
continuous on the set where $y \ge 0$.

	Similarly, if the Fourier transform of $f$ is supported in
$(-\infty, 0]$, then the integral for the inverse Fourier transform of
$f$ defines a holomorphic function on the lower half-plane, consisting
of the complex numbers $z = x + i \, y$ with $x, y \in {\bf R}$ and $y
< 0$.  This can also be identified with the Abel sum for the inverse
Fourier transform with $y = -t$.  As a result one obtains a continuous
function on the set of $z = x + i \, y$ with $x, y \in {\bf R}$ and $y
\le 0$.

	Let us define $A^+_t(x)$, $A^-_t(x)$ by setting $A^+_t(x)$
equal to $A_t(x)$ when $x \ge 0$ and to $0$ when $x < 0$, and setting
$A^-_t(x)$ equal to $A_t(x)$ when $x < 0$ and to $0$ when $x \ge 0$.
Thus $A^+_t(x) + A^-_t(x) = A_t(x)$ for all $x \in {\bf R}$,
$A^+_t(x)$ is supported in $[0, \infty)$, $A^-_t(x)$ is supported in
$(-\infty, 0]$, and $A^+_t(x)$, $A^-_t(x)$ are continuous except for a
jump discontinuity at $0$.  These functions are integrable and we can
define their Fourier transforms in the usual manner, and these can be
computed directly to be $(2 \pi (t + i \, \xi))^{-1}$, $(2 \pi ( t - i
\, \xi))^{-1}$, respectively.  The sum of these two functions is the
Fourier transform of $A_t(\xi)$, which is the Poisson kernel with $n =
1$.

	Now suppose that $f(x)$ is a continuous integrable function on
the real line whose Fourier transform is supported in a closed and
bounded interval $[a, b]$.  In particular, the Fourier transform of
$f$ is integrable, and $f(x)$ is equal to the inverse Fourier
transform of its Fourier transform, which is to say as the integral of
$\widehat{f}(\xi)$ times $\exp (2 \pi i \, x \, \xi)$ with respect to
$\xi$.  Notice that this implies that $f$ is bounded.  In fact the
integral of $\widehat{f}(\xi)$ times $\exp (2 \pi i \, z \, \xi)$
makes sense for all complex numbers $z$ and defines an entire
complex-analytic function of $z$.  

	One can check directly that for each complex number $z = x + i
\, y$, $x, y \in {\bf R}$, the absolute value of this function is less
than or equal to the integral of $|\widehat{f}(\xi)|$ times the
exponential of a constant times $|y|$.  To be more precise, one can
take the exponential of the maximum of $- 2 \pi \, t \, y$ with $a \le
t \le b$.  For $y$ fixed, one gets a bounded function of $x$, and then
a linear exponential rate of growth for this as a function of $y$.

	It turns out that there are some general results about complex
analytic functions which imply analogous behavior.  Namely, suppose
that $f(z)$ is an entire complex analytic function on the complex
plane, and that $\sigma$ is a positive real number such that for each
$\epsilon > 0$ there is a positive real number $A_\epsilon$ with the
property that $|f(z)|$ is less than or equal to $A_\epsilon$ times
$\exp ((\sigma + \epsilon) |z|)$ for all $z \in {\bf C}$.  Assume
further that $f(z)$ is bounded on the real line, which is to say that
there is a positive real number $B$ such that $|f(x)| \le B$ for all
$x \in {\bf R}$.  Then $|f(x + i \, y)|$ is less than or equal to $B$
times $\exp (\sigma \, |y|)$ for all $x, y \in {\bf R}$.

	As discussed on p124-5 of \cite{4} and p108-9 of \cite{5},
this can be derived from a Phragmen--Lindel\"of type of argument.
Actually, one works on the two half-planes separately, and there are
analogous results involving other cones and growth rates.  At any
rate, one gets this nice conclusion for entire holomorphic functions
which mimics the behavior of functions whose Fourier transforms have
compact support.  This is in fact a step in the proof of certain
Paley--Wiener theorems to the effect that one can go in the other
direction and say that an entire holomorphic function with suitable
properties is the Fourier transform of a function or distribution with
compact support.

	For that matter one can consider functions on the real line
with a holomorphic extension to some region in the complex plane in
connection with restrictions on the Fourier transform.  A basic
instance of this occurs with functions which extend holomorphically to
the upper or lower half-plane in connection with Fourier transforms
supported in $[0, \infty)$ or $(-\infty, 0]$, as before.  One can look
at this in terms of Cauchy integral representations, which can be
described in terms of convolutions, and analyzed further using the
fact that the Fourier transforms of $A^+_t$, $A^-_t$ correspond
exactly to Cauchy kernels for the upper and lower half-planes as
discussed previously.  Much more information on these topics can be
found in the references.

\end{document}